\font\fourteenbf=cmbx10  scaled\magstep2
\vglue 0.5truein

 {\fourteenbf On Lange's Conjecture}
\bigskip

\bigskip
Montserrat Teixidor i Bigas,

Mathematics Department, Tufts University,
Medford MA 02155, U.S.A., mteixido@tufts.edu

temporary address DPMMS, 16 Mills Lane, Cambridge CB2 1SB, England.
teixidor@dpmms.cam.ac.uk

\bigskip
Let $C$ be a generic curve of genus $g$ defined over an algebraically
closed field of characteristic zero. Denote by $U(n,d)$ the moduli space
of stable vector bundles of rank $n$ and degree $d$ on $C$. Let $E$
be an element of $U(n,d)$. Consider the set of subbundles $E'$ of $E$
of rank $n'$. Define the integer $s_{n'}(E)=n'd-nmax\{deg E'\}$.
Then, $E$ is semistable if and only if $s_{n'}(E)\ge 0$ for all $n'$.
On the other hand, a generic $E\in U(n,d)$ must satisfy 
$s_{n'}(E)\ge n'(n-n')(g-1)$ (cf [L] Satz p.448), while for any $E$, 
$s_{n'}(E)\le n'(n-n')g$ (cf [L] Satz p.448, [M-S]).

One can then stratify $U(n,d)$ according to the value of $s_{n'}$: 
define 
$$U_{s,n'}(n,d)=\{ E\in U(n,d)|s_{n'}(E)\le s \}, 0\le s\le n'(n-n')g$$
These loci are closed in $U(n,d)$ because the function 
$s_{n'}(E)$ is upper semicontinuos on $E$. It is not clear though
whether the loci are non-empty. The most natural approach
towards the study of $U_{s,n'}(n,d)$ is to consider families of 
extensions of  vector bundles  $E''$ by vector bundles $E'$. 
Here $E''$ is assumed to have rank $n-n'$ and degree $d-d'$.
If $\mu (E')<\mu (E'')$, one expects the generic such extension
to be stable. One could even restrict attention to the subloci
$\bar U_{s,n'}(n,d)$ of vector bundles that can be written as an
 extension with both $E',E''$ semistable. A moduli count shows 
that the dimension of the set of such extensions is
 $d_{s,n'}=n^2(g-1)+1+s-n'(n-n')(g-1)$.
 Hence, one expects $dim U_{(s,n')}(n,d)=
d_{s,n'}$ if $0<s\le  n'(n-n')(g-1)$ and $U_{(s,n')}(n,d)=U(n,d)$
if $s\ge n'(n-n')(g-1)$. 
Lange (cf [L]) conjectured that $\bar U_{s,n'}(n,d)$ should be 
non-empty for $0<s\le n'(n-n')(g-1)$. 
 The conjecture had been proved so far for vector bundles of rank
two (cf[L,N]) and for vector bundles of any rank if 
$0<s<min(n',n-n')(g-1)$ and in a few very special cases (cf [B,B,R]).
 The method used in these papers consists in showing that the set of 
unstable extensions fill a space of dimension smaller than $d_{(s,n')}$.
After this work was completed, we received a preprint from 
H.Lange and L.Bambrila-Paz proving results similar to ours.
 
We take a different approach which continues our study in [T3].
 We work with a curve generic in the sense of moduli.
We construct explicitly
stable vector bundles of degree $d$ and rank $n$ that have subbundles
of rank $n'$ and degree $d'$ in the suitable range for 
a special type of reducible curve. We then show that this is enough
in order to prove the result for the generic curve.
Our method  also allows us to clarify what happens in the so far grey area
$n'(n-n')(g-1)\le s\le n'(n-n')g$. The situation is as follows: subbundles of 
rank $n'$ and degree $d'$ exist on the generic $E$ as soon as 
$s\ge n'(n-n')(g-1), s\equiv n'd (n)$. Our result is
\bigskip

\proclaim (0,1) Theorem. Let C be a generic curve of genus $g\ge 2$ 
defined over a field of characteristic zero. Fix integers
$n,d,0<n'<n,s\equiv n'd(n)$. Define $d'=(n'd-s)/n$.
\smallskip
i)If $0\le s\le n'(n-n')(g-1)$, then $\bar U_{s,n'}(n,d)$ is non-empty and
there are elements $E\in \bar U_{s,n'}(n,d)$  that have only a finite
 number of subbundles of rank $n'$ and maximal degree.
\smallskip
ii)Assume $n'(n-n')(g-1)\le s$. Let $E$ be a generic vector bundle
 on $C$ and $Q$ a generic point on $C$.
Then $E$ can be written as an extension 
$0\rightarrow E'\rightarrow E\rightarrow E''\rightarrow 0, rkE'=n', deg E'=d'$
and both $E',E''$ semistable.
 The set $A=A(n,,d,,E)$ of subbundles of $E$ of rank 
$n'$ and degree $d'$ is non-empty of dimension $s-n'(n-n')(g-1)$ and the map
$$\matrix {\pi _Q:&A&\rightarrow &{\bf G}(n',E_Q)\cr
 &E'&\rightarrow& E'_Q\cr}$$
is either onto or has generic finite fibers.
                                                             
\bigskip
Note: for elliptic curves, the result is essentially correct (see (2,8))
except that stability must be replaced by semistability. This is to
 be expected as a ``generic'' vector bundle on an elliptic curve is 
only semistable and not stable if rank and degree are not coprime 
(see [Tu] for an alternative approach to moduli spaces of vector 
bundles on elliptic curves).

The paper is structured as follows: in section 1 we gather a few technical
results. These will allow us to reduce the proof of the Theorem for the 
generic curve to the proof in the case of a special reducible curve.
Specifically, we are going to deal with the curve that can be
obtained by gluing at a point a generic curve of genus  one
less with an elliptic curve. 
 Section 2 proves (a strenghtening of)
the Theorem for the case of an elliptic curve. In addition to providing 
the first step of the induction, it will give us the working tool
to allow us to treat the tail of our reducible curve. This section 
makes heavy use of Atiyah's classical paper [A]. We need here 
characteristic zero (see the proof of (2,8)) although the results are 
likely to hold in any characteristic.
The final section will prove the Theorem for the reducible curve 
under the assumption that it holds for the generic curve of smaller genus.

Acknowledgements : This work has been carried out during a visit
to the Pure Mathematics Department at the University of Cambridge.
The author is a member of the Europroj group
``Vector Bundles on Curves''. The interest in this question 
was spurred by a talk given by Barbara Russo at the Europroj
meeting 96. Thanks are also due to Peter Newstead for pointing
out an incorrect definition of $s_{n'}(E)$.

\bigskip

{\bf 1 Reduction to the case of a reducible curve.}
\bigskip

We mainly use standard notation. We gather below a few definitions, 
including some that have already been introduced above.

{\bf (1,1) Notation.} As above $U(n,d)$ denotes the moduli space of stable
vector bundles on a curve. When the curve is reducible, 
stability is understood with respect to a polarisation
$(a_i)$ (cf.[S]).
 $U_{s,n'}(n,d)$ denotes the set of stable
vector bundles that have a subbundle of rank $n'$ and degree $d'$ 
with $n'd-nd'=s$. When there is some danger of confusion, we shall 
add reference to the curve as  $U_{s,n'}(n,d)(C)$.  

$\bar U_{s,n'}(n,d)$ denotes the set of stable
vector bundles that have a stable subbundle of rank $n'$ and degree $d'$ 
with $n'd-nd'=s$ and the quotient is also stable.

$A(n',d',E)$ will denote the set of subbundles of rank $n'$ and degree
$d'$ of a vector bundle $E$ on the curve $C$. When $C$ is reducible, 
we shall often replace the degree by a multidegree (or degree of the 
restricition to the different components). If the curve is $C_k, C^k$ instead
of $C$, we shall write the sub or superindex $k$ in $A$ as well.

${\bf G}(n',V)$ denotes the Grassman variety of linear subspaces of 
dimension $n'$ of a vector space $V$ of dimension $n$.

$\pi_P:A(n',d',E)\rightarrow {\bf G}(n',E_P)$ denotes the (rational)
map that sends a subbundle to its fiber at a (generic) point $P\in C$

\bigskip

\proclaim (1,1) Lemma. Let $C$ be a projective connected reduced
curve having only nodes as singularities. Let 
$0\rightarrow E'\rightarrow E\rightarrow E''\rightarrow 0$
be an exact sequence of vector bundles. Assume that the restriction 
to each component of $C$ of each of the vector bundles above 
is semistable and that $\mu (E')< \mu (E)$. Then $h^0(E''^*
\otimes E')=0$.

Proof: A non-zero section of $E''^*\otimes E'$ is equivalent to a 
non-zero map $f:E''\rightarrow E'$. As the restriction of $E$
to each component $C_i$ is semistable, $\mu (E'_{|C_i})\le \mu (E''_{|C_i})$
Moreover, as $\mu (E')<\mu (E)$, there is at least one component 
$C_{i_0}$ such that  $\mu (E'_{|C_{i_0}})< \mu (E''_{|C_{i_0}})$.
By the semistability of the bundles involved, $f_{|C_{i_0}}=0$.
As $C$ is connected, given any component $C_{i_1}$, there exists a chain 
$C_1,...,C_t$ such that $C_1=C_{i_0}, C_t=C_{i_1}$ and $C_i$ intersects
$C_{i-1},i=2...t$. We prove by induction on $t$ that $f_{|C_t}\equiv 0$.
If $t=1$, we know this already. Assume $f_{|C_{t-1}}=0$. Then the 
map $f_{|C_t}$ vanishes at a point $P$ in the intersection of 
$C_{t-1}$ and $C_t$. Hence it gives rise to a map
 $\bar f_{|C_t}:E''(P)_{|C_t}\rightarrow E'_{|C_t}$.
As $\mu (E''(P)_{|C_t})>\mu (E'_{|C_t})$ the semistability of the bundles 
implies that $\bar f_{|C_t}=0$. Hence $ f_{|C_t}=0$. This concludes the induction and hence the proof of the Lemma.

\bigskip

The following proposition and section 3 deal with vector bundles
on reducible nodal curves. There is a notion of stability
for torsion-free sheaves on such a curve introduced by 
Seshadri. We refer the reader to [S] Ch.VII-VIII, where he 
constructs a moduli spaces for such objects.
Here we recall how to determine a vector bundle on such 
a curve. 

Let $C$ be a nodal curve with irreducible components $C_i$
and nodes $P_j$. The restricion of a vector bundle $E$
to each $C_i$ is a vector bundle $E_i$. If $P_j$ is in the intersection 
of $C_{i_1}$ and $C_{i_2}$, there is a ``gluing'' 
isomorphism of the  vector spaces fiber at $P_j$:
$E_{i_1,P_j}\cong E_{i_2,P_j}$. Conversely, a set of such data determines
a vector bundle. A torsion-free sheaf is determined by giving isomorphisms
of suitable quotients of the fibers
 $E_{i_1,P_j}/V_{i_1,j}\cong E_{i_2,P_j}/V_{i_2,j}$.
\bigskip

\proclaim (1,3) Proposition. Let $C^0$ be a tree-like curve
(or curve of compact type). Let
 $$0\rightarrow E'^0 \rightarrow E^0 \rightarrow E''^0\rightarrow 0$$ 
be an exact sequence of vector bundles on $C^0$. Assume that the 
restriction of each vector bundle to each component of $C^0$ is semistable
and stable on at least one of them. Assume also that the restriction 
of the vector bundles to any destabilising rational component is 
trivial (i.e. of the form ${\cal O}^n$).
 Denote by $n,d,n',d'$ the rank and degree of 
$E^0, E'^0$ respectively, $s=n'd-nd'$.\smallskip
i) Assume $0< n'd-nd'\le n'(n-n')(g-1)$,
and $E^0$ has only a finite number of subbundles of rank $n'$
 and degree $d'$. Then, $\bar U_{s,n'}(n,d)(C)$ is non-empty
for the generic curve $C$ of genus $g$ 
and the generic element has only a finite number of subbundles
of rank $n'$ and degree $d'$.\smallskip
ii) Assume $n'(n-n')(g-1)\le s$ and
 $dim A^0(n',d',E^0)=n'd-nd'-n'(n-n')(g-1)$.
 Then the generic element $E$ in $U(n,d)(C)$ has a semistable
subbundle of rank $n'$ and degree $d'$ with semistable quotient and
 $dim A(n',d',E)=n'd-nd'-n'(n-n')(g-1)$.

Proof: This proof is inspired in  Lange's construction of $\bar U_{s,n'}$ (cf.
[L] section 4) adapted to the case of a family and is essentially contained in 
[T3]. We include it here for the sake of completeness.

 Claim: The set of vector bundles that are trivial on rational components
 and satisfy the stability conditions stated above can be identified
 to an open set of the moduli space of torsion -free 
sheaves on the stable model $\tilde C$ of $C^0$:
Consider a vector bundle on $C^0$, that is trivial on a rational component,
There is a canonical isomorphism between the fibers at the two
 nodes $P_1,P_2$ on the rational component induded by the pair of 
isomorphisms $H^0({\bf P^1}, {\cal O}^n)\rightarrow {\cal O}^n_{P_i} i=1,2$.
 This induces a natural gluing isomorphism 
for the vector bundle on the curve with this component blown down
. Proceed in this way until all components are blown down.
 Conversely, given a vector bundle (say $\tilde E$) on
 $\tilde C$, one can build a vector 
bundle on $C^0$ by taking trivial bundles ${\cal O}^n$ on the additional
 rational components. The gluings at the nodes are canonically determined 
by those of $\tilde E$. 

Conversely, we can replace $C^0$ by a curve $C'^0$ with chains of rational 
components added between some of the nodes of $C^0$
(i.e. a semistable curve with same stable model as $C^0$ 
a little bit more unstable than $C^0$).

Consider  a family ${\cal C}\rightarrow S$ parametrised say 
by the spectrum of a discrete valuation ring whose central fiber 
is $C^0$ and whose generic fiber is a non-singular curve $C$ of genus
$g$. From the reasoning above, adding rational components to 
$C^0$ does not alter the hypothesis. Hence, we can make base changes and 
normalise if needed.

We now need the existence of moduli spaces parametrising vector bundles 
 of the type $E'^0, E^0, E''^0$ on $\cal C$. In fact we need suitable 
coverings of these moduli spaces with Poincare families.
Denote by $( C_i)_{i\in  I}$ the irreducible components of $ C^0$.
Choose a polarisation $(a_i)_{i\in  I}$ for $ C^0$.
 This is the choice 
of rational numbers (``weights'') $a_i>0, \sum a_i=1$.
Consider moduli spaces of torsion-free sheaves on ${\cal C}
\rightarrow S$  of ranks $n',n,n-n'$ and degrees $d',d,d-d'$
such that the restriction to the 
generic fiber is semistable and the restriction to 
the special fiber is $(a_i)$-semistable. 
These can be shown to exist as in Seshadri's proof for a 
single curve by replacing points by sections
of the family. Alternatively, one can use the existence of 
a moduli space for torsion-free sheaves over the moduli space 
$\bar {\cal M}_g$ of stable curves constructed by Pandharipande([P]).
Our curves are only semistable and not necessarily stable 
but our vector bundles are trivial on rational components. 
By blowing down these rational components, we can identify the open set 
of these bundles with an open set of the moduli space of vector bundles
on the stable model of the curve as explained above. Take then in the
standard way suitable coverings over  which Poincare bundles exist.
Denote by ${\cal U}', {\cal U}, {\cal U}''$ these moduli spaces, by
 ${\cal M'}, {\cal M}, {\cal M''}$ their coverings and by 
${\cal E'}, {\cal E}, {\cal E''}$ the corresponding Poincare bundles.
The given  vector bundles $E'^0,E^0,E''^0$ are not necessarily
$(a_i)$-stable and therefore may not be represented in these families
. In order to deal with this situation, 
we need to recall some results about stability of vector bundles 
on tree-like curves  (see [T1], [T2]).

Assume a vector bundle on $C^0$ is $(a_i)$-semistable. One can choose
an ordering of the components of $C^0$ in such a way that the 
complement  of the components chosen beforehand is 
connected. Then the degree 
on each component lies on intervals of lenght $n$ in the real line
determined by d and the degree in the previous components.
We shall call such a multidegree an $(a_i)$-allowable multidegree.
The converse is not true in general: a vector bundle with an $(a_i)$-
allowable multidegree is not necessarily stable. However if 
$E$ is a vector bundle with an $(a_i)$-allowable multidegree, 
 the restriction to each 
component is semistable and at least one of them is stable, 
then $E$ is stable.
One can choose a multidegree on $C^0$ such that for every 
divisor $ D^0$ on $ C^0$ with this multidegree, $ E^0( D^0)$
has its multidegree in these intervals. Everything remains
valid replacing $E^0$ by say $E'^0,E''^0$ or replacing 
$C^0$ by $\tilde C$.

As $C^0$ is of compact type, one can choose a divisor
${\cal D}$ on ${\cal C}$ with support on $C^0$ whose restriction to 
$C^0$ has a given preassigned degree. By the normality 
of ${\cal C}$, ${\cal D}$ corresponds to a line bundle
${\cal L}$.  Similarly, there are 
line bundles ${\cal L'}, {\cal L''}$ such that 
$E'^0\otimes {\cal L'}_{|C^0},
 E''^0\otimes {\cal L''}_{|C^0}$ 
are $(a_i)$-semistable. Then ${\cal M'},{\cal M}, {\cal M''}$
 with Poincare bundles ${\cal E'}\otimes {\cal L'}^{-1},
{\cal E}\otimes {\cal L}^{-1}, {\cal E''}\otimes ({\cal L''})^{-1}$ 
parametrise families of vector bundles 
containing $E'^0,E^0, E''^0$. Moreover, there are open sets  of 
${\cal M'},{\cal M}, {\cal M''}$ such that the restriction of  the 
corresponding vector bundle to each component of $C^0$ is semistable.
Denote by $q, p',p''$ the projections of 
${\cal C}\times _S{\cal M}'\times _S {\cal M''}$
on the three factors. Consider
$${\cal F}=R^1(p'\times p'')_*[(q\times p'')^*({\cal E''} \otimes 
({\cal L''})^{-1})^*\otimes (q\times p')^*({\cal E'}
 \otimes ({\cal L'})^{-1})].$$
From (1,2) and Riemann-Roch's Theorem, ${\cal F}$ is a vector bundle
of rank $n'(n-n')(g-1)+n'd-nd'$ over ${\cal M'}\times {\cal M''}$. 
 Then, ${\bf P}({\cal F})$ parametrises extensions of vector bundles
and there is a universal family (cf [L2])
$$0\rightarrow (id\times \pi)^*({\cal E'}\otimes {\cal L}'^{-1})
\rightarrow \tilde {\cal E}\rightarrow 
(id\times \pi)^*({\cal E''}\otimes {\cal L}''^{-1})\rightarrow 0 $$
where $\pi$ denote pull-backs to${\bf P}({\cal F})$. 
Using the $(a_i)$-allowability of the multidegrees of the endterms
of this sequence, one can check the $(a_i)$-allowability of the 
central term. As (semi)-stability is an open condition and
the restriction of $E^0$ to each component of $C^0$ is 
semistable, the same is true for the generic extension. 
Hence the generic extension is $(a_i)$-stable. There is
then a rational map ${\bf P}({\cal F})\rightarrow U$. 
If $s\le n'(n-n')(g-1)$, our assumption implies that the fibers
of the map are finite. Hence the image has dimension 
$n'^2(g-1)+1+(n-n')^2(g-1)+1+n'(n-n')(g-1)+n'd-nd'-1=n^2(g-1)+1
+s-n'(n-n')(g-1)$ as claimed. If  $s\le n'(n-n')(g-1)$,
 our assumption implies that the fibers
of the map have dimension at most $s-n'(n-n')(g-1)$. Hence,
 the map is generically onto as claimed.
 
{\bf Remark}. The condition that the restriction of the 
vector bundles to one of the components of $C^0$ is stable 
could be weakened. All we need in fact is that the 
vector bundle one obtains on the whole curve is stable by some polarisation
(up to tensoring each restriction with a line bundle).

\beginsection 2 Subbundles of vector bundles on elliptic curves.

In this section C will always be an elliptic curve. A vector bundle is
said to be indecomposable if it is not the direct sum of smaller rank 
subbundles.  Up to the choice 
of a point $P$ in $C$, there is a canonically determined indecomposable
vector bundle $E(n,d)$ of rank $n$ and degree $d$ on 
$C$ (cf. [A]). When $d=0$, $E(n,0)$ is the only indecomposable 
vector bundle that has sections. Every indecomposable vector bundle on $C$
is of the form $E(n,d)\otimes L$ for some line bundle $L$
of degree zero. Hence $Pic^0(C)\cong C$ can be taken as
a parameter space for the family of indecomposable vector 
bundles on $C$. Denote by ${\cal L}$ the Poincare bundle of degree
zero on $C\times C\cong C\times Pic^0(C)$. Then 
${\cal E}(n,d)=p_2^*(E(n,d))\otimes {\cal L}$ is a Poincare bundle 
on $C\times C$
(where $p_2$ denotes projection onto the second factor).

 Every vector bundle on $C$ is a direct sum
 of indecomposable vector bundles.
If we fix a decomposition type, say
 $$ (2,1)d^j, n^j>0,j=1...k, \sum_j d^j=d, \sum_jn^j=n,$$
one can take as parameter space for vector bundles with this
type of decomposition the $k^{th}$ product of the curve $C$. 
The corresponding Poincare bundle is 
$$(2,2){\cal E}((n^j,d^j))=\oplus [p_{k+1}^*(E((n^j,d^j)))\otimes
 p^*_{j,k+1}({\cal L})].$$
Here again $p$ denote projections and the subindices denote
the factors we are projecting onto.

\bigskip

\proclaim (2,3) Lemma.  Let $E,E'$ be 
indecomposable vector bundles on $C$ of ranks $n,n'$ and degrees
$d,d'$ respectively. Then, $E'^*\otimes E$ is a direct sum of 
indecomposable vector bundles all of them with the same slope.

\bigskip
Proof: From [A] p.434 Corollary iv), $E(n',d')^*=E(n',-d')$. 
If $(n,d)=(n',d')=1$, the result is then Lemma 33 of [A]. 
If $(n,d)=h, (n',d')=h', n=\bar nh, n'=\bar n'h'$ then by [A]
Lemma 24, $E(n,d)=E(h,0)\otimes E(n_1,d_1),
E(n',d')=E(h',0)\otimes E(n'_1,d'_1)$. From [A] Lemma 21, 
$E(h,0)\otimes E(h',0)$ is a direct sum of indecomposable vector bundles 
of degree zero. Then, the result follows using [A] Lemma 33
aplied to $E(\bar n,\bar d)\otimes E(\bar n',\bar d')$ and [A] Lemma 24
again.

\bigskip
\proclaim (2,4)Lemma. Let $E$ be an 
indecomposable vector bundle of rank $n$ and degree $d$.
Fix integers $n'>0,d'$. Consider all possible decomposition types
for a vector bundle $E'$ of rank $n'$ and degree $d'$: 
$ d_i^j, n_i^j>0,j=1...k_i, \sum_j d_i^j=d', \sum_jn_i^j=n'$.
 Let $X_i$ be the set of points 
in $C\times...\times C$ representing vector bundles $E'_i$ such that 
there is a map from $E'$ to $E$ that is non-zero on each indecomposable 
summand $E'^j_i$ of $E'_i$. There is only a finite number of $X_i$.
Define ${\bf Hom_i}=p_{1...k_i*}[{\cal E}((n_i^j,d_i^j))^*\otimes 
p_{k_{i+1}}^*E]_{|X_i}$. Then, ${\bf Hom_i}$ is a vector bundle
$$(2,5)dim X_i+rk {\bf Hom_i}\le n'd-nd'+h^0(E'^*_i\otimes E'_i)$$
where $E'_i$ is a generic point of $X_i$. 
 Define $(n',d')=h',
n'=\bar n'h', d'=\bar d'h'$ and $\bar X$ as the component
corresponding to the decomposition of $E'$ into direct sum of 
$h'$ subbundles of degree $\bar d'$ and rank $\bar n'$. Then on $\bar X$
(2,5) is an equality except may be in the case when $n'd-nd'=0, h'>0$ 
. Equality in (2,5) implies that all the sumands of $E'$ have the same slope.
If $(n,d)=1$ or $n'd-nd'>0$, equality in (2,5) occurs
only on $\bar X$.

\bigskip
Proof: We can write $E=E(n,d)\otimes L$. 
A vector bundle of degree $d_i^j$ and rank $n_i^j$ is of the form 
$E(n_i^j,d_i^j)\otimes L_i^j$ for a degree zero line bundle $L_i^j$.
A generic point $E'_i$ of $X_i$ is of the form
 $E'_i=\oplus_j(E(n_i^j,d_i^j)\otimes L_i^j)$.  
From (2,3), $(E(n_i^j,d_i^j)\otimes L_i^j)^*\otimes E(n,d)\otimes L$
 is a direct sum of indecomposable 
vector bundles all of them of slope $(n_i^jd-nd_i^j)/n_i^jn$.
If $n_i^jd-nd_i^j<0$, then
 $(E(n_i^j,d_i^j)\otimes L_i^j)^*\otimes (E(n,d)\otimes L)$ has no nontrivial
sections (use for instance [A] p.434 Corollary iv)and [A] lemma 15).
Hence  $X_i$ is empty.
Therefore, we need only consider a finite number of decomposition types
corresponding to $n_i^jd-nd_i^j\ge 0$ for all possible choices of $j$.

From (2,3) above and [A] Lemma 15, if $n_i^jd-nd_i^j>0$
$$h^0([E(n_i^j,d_i^j)\otimes L_i^j]^*\otimes [E(n,d)\otimes L]
=n_i^jd-nd_i^j.$$ 

Asume now $n_i^jd-nd_i^j=0$. Write
 $$(n,d)=h, (n_i^j,d_i^j)=h_i^j, n=h\bar n, d=h \bar d, 
n_i^j=h_i^j\bar n_i^j, d_i^j=h_i^j\bar d_i^j .$$
Then
 $$E(n,d)=E(h,0)\otimes E(\bar n,\bar d), E(n_i^j,d_i^j)=
E(\bar n_i^j,\bar d_i^j)\otimes E(h_i^j,0)$$
As $0=n_i^jd-nd_i^j=hh_i^j(\bar n_i^j \bar d-\bar n \bar d_i^j)$ and 
$(\bar n,\bar d)=1$, it follows that $\bar n|\bar n_i^j$. 
Similarly, $\bar n_i^j|\bar n$. Hence, $\bar n_i^j=\bar n, \bar d_i^j=
\bar d$. Then 
$$(E(n_i^j,d_i^j)\otimes L_i^j)^*\otimes (E(n,d)\otimes L)=E(h_i^j,0)\otimes
E(h,0)\otimes E(\bar n, \bar d)\otimes E(\bar n, \bar d)^*\otimes L
\otimes (L_i^j)^{-1}.$$

 Using [A] Lemma 22, $E(\bar n, \bar d)\otimes
E(\bar n, \bar d)^*=\oplus L_k$ where the $L_k$ are the line bundles
of order  dividing $\bar n$.
Hence, from [A] Lemma 21,
 $$E(n_i^j, d_i^j)^*\otimes E(n,d)=
[E(|h-h_i^j|+1,0)\oplus E(|h-h_i^j |+3,0)\oplus...\oplus E(h+h_i^j-1,0)]
\otimes 
\oplus(L_kL(L_i^j)^{-1}$$

 Hence the space of maps from 
$E(n_i^j,d_i^j)$ to $E(n,d)$ is non-empty if and only if $L_i^j=LL_k$ for some 
k. This means that there is only a finite number of choices for 
$L_i^j$. For this finite number of choices , the dimension of
the space of maps from $E(n_i^j,d_i^j)\otimes L_i^j $ to $E(n,d)$ is
 $min(h,h_i^j)$. 

Therefore, $h^0(E'^*_i\otimes E)$ is constant for $E'_i$ on a fixed 
$X_i$. Hence, ${\bf Hom_i}$ is a vector bundle of rank $h^0(E'^*_i\otimes E)$.

 Write $\epsilon _i^j=0$ if 
$n_i^jd-nd_i^j>0$, $\epsilon _i^j=1$ if $n_i^jd-nd_i^j=0$.
Then 
$$dim X_i+rk {\bf Hom_i}=dim X_i+ \sum (n_i^jd-nd_i^j)+\sum \epsilon ^j_i
 min(h,h_i^j)=$$
$$=(k_i-\sum \epsilon _i^j)+ n'd-nd'+\sum \epsilon_i^j min(h,h_i^j).$$

We now compute the dimension of the space of automorphisms of $E'_i$.
From the description of $E'_i$, 
$$E'^*_i\otimes E'_i=\oplus_j[E(n_i^j,d_i^j)^*\otimes E(n_i^j,d_i^j)]\oplus
\oplus_{j_1\not= j_2}[(E(n_i^{j_1},d_i^{j_1})\otimes L_i^{j_1})^*
\otimes (E(n_i^{j_2},d_i^{j_2})\otimes L_i^{j_2}]$$ 
 
Write $(n_i^j,d_i^j)=h_i^j$ and $n_i^j=h_i^j\bar n_i^j, d_i^j=h_i^j\bar d_i^j$. Then $E(n_i^j,d_i^j)=E(\bar n_i^j,\bar d_i^j)\otimes E(h_i^j,0)$ 
(cf[A], Lemma 24). From [A] Lemma 22, $E(\bar n_i^j,\bar d_i^j)
\otimes E(\bar n_i^j,\bar d_i^j)^*=\oplus L_i^{jk}$
 where the $L_i^{jk}$ are the linebundles on $C$ of order 
dividing $\bar n_i^j$. From [A] Lemma 21,
 $$E(h_i^j,0)\otimes E(h_i^j,0)^*=E(1,0)\oplus E(3,0)\oplus ...\oplus
 E(2h_i^j-1,0).$$
 Hence, $h^0(End E(n_i^j, d_i^j))=h_i^j$ . 

If $ n_i^{j_2}d_i^{j_1}-n_i^{j_1}d_i^{j_2}>0$, then, using (2,3), 
$$h^0((E(n_i^{j_2},d_i^{j_2})\otimes L_i^{j_2})^*
\otimes (E(n_i^{j_1},d_i^{j_1})\otimes L_i^{j_1}))=
n_i^{j_2}d_i^{j_1}-n_i^{j_1}d_i^{j_2}$$

If $n_i^{j_2}d_i^{j_1}-n_i^{j_1}d_i^{j_2}=0$, then reasoning as above,
there is only a finite number of choices of $L_i^{j_2}$ for which 
the space of maps to the component corresponding to $j_1$ is
 non-zero. For a generic point of $\bar X$ the $L_i^j$ are generic
except when $n'd-nd'=0$. So, for a generic point of $\bar X$, 
 all these spaces are zero except may be in the case when $n'd-nd'=0$.
Hence,
$$h^0(E'^*_i\otimes E'_i)\le \sum_j h_i^j+\sum _{j_1<j_2}
|n_i^{j_2}d_i^{j_1}-n_i^{j_1}d_i^{j_2}|.$$
with equality for the generic $E$ in $\bar X$ if $n'd-nd'>0$.
As $h_i^j\ge 1$ for all j, $dim X_i+rk {\bf Hom_i}
\le n'd-nd'+h^0(E'^*_i\otimes E'_i)$ 
as stated. Moreover, if equality occurs, 
then  $n_i^{j_2}d_i^{j_1}-n_i^{j_1}d_i^{j_2}=0$ for all pairs
$j_1,j_2$. Also, $h_i^j=1$ if $\epsilon_i^j=0$. The first condition
implies that all summands of $E'_i$ have the same slope $d'/n'$.
 If  $d'/n'\not= d/n$,   then  $h_i^{j_1}=h_i^{j_2}=1$. This 
together with the equality above
implies $d_i^{j_1}=d_i^{j_2}, n_i^{j_1}=n_i^{j_2}$. Hence,
 we are dealing with a point of $\bar X$. Similarly, if $(d,n)=1$,
then $min(h,h_i^j)=h=1$ and equality implies that $h_i^j=1$. 
We conclude as before that we are dealing with a point of $\bar A$.

\bigskip
\proclaim (2,5) Lemma. Let $d,n,d', n', 0<n'<n$
 be integers, $h=(n,d), d=\bar dh,n=\bar nh$. Let 
$E=\oplus E_i$ be a vector bundle written as the direct sum of 
$h$ generic indecomposable vector bundles of rank $\bar n$ and degree
$\bar d$. There is a variety $A=A(n',d',E)$ parametrizing the set of subbundles
 of $E$ of rank $n'$ and degree $d'$. The variety  $A$ is a finite union 
of  components all of them of dimension
at most $n'd-nd'$ (and is empty if $n'd-nd'<0$). The component
 $\bar A$ corresponding to an $E'$
in $\bar X$ (as defined in (2,4))
has precisely this dimension.
\bigskip

Proof: Consider the possible splitting type of a bundle $E'$ of 
rank $n'$ and degree $d'$. As an immersion cannot vanish on any
direct summand of $E'$, only a finite number 
of these decomposition types are possible (cf.(2,4)).
With the notations of (2,4), $E'$ corresponds to a point in 
some $X_i$. Let $X$ be the (disjoint) union of the components 
$X_i$ in (2,4) such that the generic map of a point in $X_i$ to $E$ is an 
immersion. Let ${\cal G}$ be the Grassmanian bundle of dim $n'$ 
subspaces of the fibers of $E$. There is then a rational 
map from ${\bf P}(\cup {\bf Hom_i})$ to ${\cal G}$. Denote by 
$A$ its image. Lange proved that $A$ can be identified to 
the image subvariety in ${\cal G}$ (cf.[L], proof of Lemma
4.2). If two maps from the same sheaf $E'$ to $E$ differ 
in an automorphism of $E'$, they give rise to the same subbundle
of $E$. Conversely, if two immersions of vector bundles have the same subbundle
as image, then, the two subbundles are isomorphic and the immersions differ
in an automorphism of the subbundle.
 Hence, from (2,4), $dim A\le n'd-nd'$ with equality 
only on the component $\bar A$ corresponding to $\bar X$ if $n'd-nd'>0$.

When $\bar n'\bar d-\bar d' \bar n=0$, then $\bar n=\bar n', \bar d=\bar d'$.
It is then easy to see that $E'$ must be the direct sum of $h'$ of the
 summands appearing in the decomposition of $E$. Hence,the result is obvious. 

It remains to show that $\bar A$ (and hence $A$) is non-empty by proving
 that the generic map from a point of $\bar X$ to $E$ is injective. This is
clear if $n'd-nd'>0$.

Assume now $n'd-nd'>0$ and the generic map from a point $E'$ in $\bar X$
to $E$ is not injective. Then the kernel $K$ is a subsheaf of $E$. As $E$
is torsion-free, $K$ is a vectorbundle of rank
$n''<n'$ and degree $d''$. The map $f:E'\rightarrow E$ can be factored 
through the canonical projection $E'\rightarrow E'/K$ and $f$ is determined 
by $K$ and the immersion $E'/K\rightarrow E$. As $E$ is torsion-free,
so is $E'/K$ and hence $E'/K$ is a vector bundle. Then, 
$$n'd-nd'=dim\bar A(n',d',E)\le dim A(n'',d'',E')+dim A(n'-n'', d'-d'',E)\le$$
$$\le n''d'-n'd''+(n'-n'')d-n(d'-d'').$$
Hence,
 $$(2,5,1)(d-d')/(n-n')\le d''/n''.$$
 On the other hand, by the semistability
of $E'$ (see for example [Tu] Appendix A),
 $$(2,5,2)d''/n''\le d'/n'$$ 
The hypothesis $n'd-nd'>0$ is equivalent to 
$$(2,5,3)d'/n'<(d-d')/(n-n').$$

 As (2,5,1),(2,5,2),(2,5,3) are incompatible, we obtain a contradiction
and the result is proved.

\bigskip

\proclaim (2,6) Lemma. With the notations of (2,5), let
 $E$ be a direct sum of vector bundles all
of the same slope. Let $E'$ be a generic point in $\bar X$ and 
$f:E\rightarrow E'$ a generic immersion. Define $E''$ by the exact
sequence
$$0\rightarrow E'\rightarrow E\rightarrow E''\rightarrow 0$$
Then $E''$ is a direct sum of indecomposable bundles all of the same slope.

Proof: There is a one to one correspondence between immersions of subbundles
$E'$ of $E$ of rank $n'$ and degree $d'$ and quotients of $E$ of rank $n-n'$
and degree $d-d'$. Hence, from (2,4), the set of such quotients has 
dimension $n'd-nd'$.

We show first that $E''$ is a vector bundle. If this were not the case, 
the immersion $E'\rightarrow E$ would factor through an immersion
$E'\rightarrow \tilde E' \rightarrow E$ with $\tilde E'/E'$
a torsion sheaf. Write $\tilde d'=deg \tilde E'=d'+k, k>0$.
There is then an exact sequence $0\rightarrow E'\rightarrow
\tilde E' \rightarrow {\bf C}^k\rightarrow 0$, where 
${\bf C}^k$ is a skyscraper sheaf of lenght $k$.
The dimension of the set of such extensions is at most $n'k$.
The dimension of the set of maps $\tilde E'\rightarrow E$ is at most 
(from (2,5)) $n'd-n\tilde d'$. Hence we find a bound
$n'd-nd'=dim\bar A(n',d',E)\le n'k+n'd-n(d'+k)=n'd-nd'-k(n-n')$. 
As $n\ge n'$, this implies $k=0$ and so $\tilde E'=E'$.

Dualizing the exact sequence defining $E''$, we obtain a one to one 
correspondence between the set of quotients $E''$ and subbundles 
of $E^*$ of rank $n-n'$ and degree $-(d-d')$. From (2,4), the
dimension of this set is at most
$(n-n')(-d)-n(-(d-d'))=nd'-n'd$ and equality implies that $E''$ is of
the stated type. This finishes the proof.

\bigskip

\proclaim (2,7) Lemma. Let $E$ be an indecomposable vector bundle with slope
$\mu \ge 1$. Then $E$ is generically generated by global sections.

Proof: We use induction on the rank $n$. For n=1, every line bundle
of strictly positive degree on an elliptic curve has sections. So, 
the result follows. Assume the result for vector bundles of rank at most
$n-1$. Let $E$ be an indecomposable vector bundle of rank $n$ and 
degree $d\ge n$. Write $d=kn+\epsilon ,0\le \epsilon \le n-1, k\ge 1$.
Let $L_k$ be a line bundle of degree $k$. If $\epsilon >0$,
 $h^0(E\otimes L_k^{-1})>0$. If $\epsilon =0$, there is one $L_k$ such that
 $h^0(E\otimes L_k^{-1})>0$. On the other hand, for any line bundle $L$
of degree at least $k+1$, $E\otimes L^{-1}$ has negative degree and 
hence no sections. There is then an exact sequence 
$$0\rightarrow L_k\rightarrow E\rightarrow E''\rightarrow 0$$

Here $E''$ has no torsion (otherwise $E$ would have a line subbundle of
 degree higher than $k$). From (2,6), $E'$ is a direct sum of vector bundles 
all of the same slope $(d-k)/(n-1)=k(n-1)+\epsilon /(n-1)\ge k$. 
By induction assumption, $E''$ is generically generated by global sections.
As $h^0(E)=h^0(L_k)+h^0(E'')$, the result follows.

\bigskip

\proclaim (2,8) Proposition. Let $E, A,\bar A$ be as in (2,5).
 Let $P, Q$ be  generic points of $C$. Denote by ${\bf G}(n',E_P)$
the Grassman variety of $n'$-subvector spaces of the fiber of 
$E$ at $P$.  Consider the rational map 
$$\matrix{\pi _P:& A&\rightarrow&{\bf G}(n',E_P)\cr
 &F&\rightarrow&F_P\cr}.$$
Define $\pi _Q$ similarly.
\smallskip
i)If $0\le n'd-nd'< n'(n-n')$ the image by $\pi _P$ of each component of $A$
has dimension at most $n'd-nd'$ and one component with image of 
this dimension and finite fibers does exist. 
\smallskip
ii)If $n'(n-n')\le n'd-nd'$, then $\pi_P$ is surjective and so the
 fibers of $\pi_P$
have dimension $n'd-nd'-n'(n-n')$.
\smallskip
iii)If $ n'd-nd'\le 2n'(n-n')$, then the product map
 $\pi_P\times \pi_Q $ has generic finite fibers.
\smallskip
iv)If $2n'(n-n')\le n'd-nd'$, then the product map 
$\pi_p\times \pi_Q$ is onto.
\bigskip

Proof: Assume first $0\le n'd-nd' \le n'(n-n')$.
 From (2,3) the dimension of every component of 
$A$ is at most $n'd-nd'$. We only need to show that the restriction
of $\pi _P$ to $\bar A$ has finite fibers. We shall prove that the differential
of $\pi _P$ is injective.

The tangent space to $A$ at a point corresponding to an immersed
vector bundle $E'$ is the space of infinitessimal deformations
of $E'$ inside $E$. These are parametrised by $H^0(E'^*\otimes E/E')$.
The tangent space to a grassmanian of subspaces of a fixed dimension of 
a space $V$ at a point corresponding to a subspace $V'$ is given 
by $Hom(V',V/V')$. In our case, $V'=E'_P, V=E_P$. We need to prove then that 
the map
$$d\pi _P: H^0(E'^*\otimes (E/E'))\rightarrow E'^*_P\otimes (E_P/E'_P)$$ 
is injective when $E'$ is generic in $\bar A$ and $P$ is generic
in $C$. 
By (2,6), $E/E'$ is a direct sum of indecomposable vector bundles all
of the same slope. Then, (2,3) implies that $E'^*\otimes E/E'$ is
a direct sum of indecomposable vector bundles all of the same slope. 
The assumption $0\le n'd-nd'<n'(n-n')$ can be written
as $0\le (d-d')n'-(n-n')d'<n'(n-n')$. Hence the slope of
each indecomposable summand of $E'^*\otimes E/E'$ is positive
 and less than one. 
From [A] Lemma 15 ii), there is an immersion 
$$H^0(E'^*\otimes E/E')\otimes {\cal O}_C \rightarrow E'^*\otimes E/E'$$
Taking the fiber at the point $P$, the map will remain injective.
This proves i).
    
Assume now $n'(n-n')\le n'd-nd'$. We need to show now that the map 
$d \pi _P$ above is surjective. From (2,3) and (2,7) $E'^*\otimes E/E'$
is generically generated by global sections. This is equivalent to ii).

If $n'd-nd'<2n'(n-n')$, from (2,3)
 $E'^*\otimes E/E' (-P-Q)$ is a direct sum of 
indecomposable vector bundles of negative slope. Hence it has no sections.
 Consider the exact sequence 
$$0\rightarrow E'^*\otimes E/E'(-P-Q)\rightarrow E'^*\otimes
E/E'\rightarrow (E'^*\otimes E/E')_{P+Q}\rightarrow 0$$
As $H^0(E'^*\otimes E/E'(-P-Q))=0$, taking homology we obtain iii).

Under the hypothesis in iv), $h^0(E'^*\otimes E/E'(-P-Q))=
 h^0(E'^*\otimes E/E')-n'(n-n')$ . Hence, $H^0(E'^*\otimes E/E')$ maps onto 
$(E'^*\otimes E/E')_{P+Q}$. This proves iv).

\bigskip

{\bf 3 Proof of the result for the reducible curve.}
\bigskip

\proclaim (3,1)Lemma. Let $n',n,0<n'<n$ be integers . Let $V_1,V_2$ be vector spaces of dimension n. Let $A_1,A_2$ be subvarieties of the Grassmannians
${\bf G}(n',V_1), {\bf G}(n',V_2)$ of dimensions $k_1,k_2$ respectively. 
Consider the set $X$ of linear isomorphisms of $V_1$ to $V_2$
that identify some element in $A_1$ with some element in $A_2$.
\smallskip
i) $dim X\le k_1+k_2+n^2-n'(n-n')$
\smallskip
ii) If $A_2$ is a finite set, then  $dim X= k_1+k_2+n^2-n'(n-n')$.
\smallskip
iii) If $dim X= k_1+k_2+n^2-n'(n-n')$, then for the generic $\varphi 
\in X$, there is only a finite number of pairs $(V'_1,V'_2)\in A_1 \times
A_2$ such that $\varphi (V'_1)=V'_2$.
\smallskip
iv) If $X=Is (V_1,V_2)$, for a generic $\varphi_0\in Is(V_1,V_2)$,
$dim\{ ( V'_1,V'_2)\in A_1 \times A_2|\varphi_0 (V'_1)=V'_2\}=
k_1+k_2-n'(n-n')$
 
\bigskip
Proof: Consider the set of triples 
$$\tilde X =\{ (\varphi, V'_1,V'_2)\in Is (V_1,V_2)\times A_1 \times A_2|
\varphi (V'_1)=V'_2\}$$
Denote by $p_0,p_1,p_2$ the projections from $Is(V_1,V_2)\times
{\bf G}(n',V_1)\times {\bf G}(n',V_2)$ to the corresponding
factors. The fiber of $p_1\times p_2 $ over any point in ${\bf G}(n',V_1)
\times {\bf G}(n',V_2)$ has dimension $n^2-n'(n-n')$. Hence, the inverse image 
$\tilde X $ of $A_1\times A_2$ by $p_1\times p_2$ has dimension 
$k_1+k_2+n^2-n'(n-n')$.
As $X=p_0(\tilde X)$, i) follows. Moreover, if there is equality
in the bound for the dimension of $X$, then the fibers of 
$p_0:\tilde X\rightarrow X$ are finite. Hence, iii) follows.
If $A_2=\{ V'_2\}$ consists of a single point, then $p_0$ induces a bijection 
$\tilde X \rightarrow X$, the inverse being given by 
$\varphi \rightarrow (\varphi, (\varphi)^{-1}(V'_2), V'_2)$. 
This proves ii).

The asumption in iv) means that $p_0:\tilde X \rightarrow Is(V_1,V_2)$
is onto. Hence the generic fiber has dimension $dim\tilde X -
dimIs(V_1,V_2)$. This is iv).
\bigskip

\proclaim (3,2) Proposition. For $g\ge 2$ let C be a curve obtained by 
gluing a generic point $P_g$ in a curve $C_g$ of genus $g$ satisfying
(0,1) with a point $P_1$ on a generic elliptic curve $C_1$. For $g=1$, 
glue two generic elliptic curves $C_1, C'_1$ at generic point $P_1,P'_1$. 
 Fix integers $n,d,n'<n,s\equiv n'd(n)$. Define $d'=(n'd-s)/n$.Then,
there exists a a polarisation $(a_i)$ of $C$ such that
\smallskip
i)If $0< s\le n'(n-n')(g-1)$, there exists a vector bundle $E$ on $C$ 
that can be written as an extension 
$$0\rightarrow E'\rightarrow E\rightarrow E''\rightarrow 0$$
with $rk E'=n', degE'=d'$, $E$ has only a finite number of
 subbundles of rank $n'$ and degree $d'$ and no subbundles of rank $n'$
and larger degree. Moreover, the restrictions of $E,E',E''$ to 
$C_1$ are semistable and the restrictions to $C_g$ are stable if 
$g>1$. If $g=1$, the restrictions to $C'_1$ are only semistable 
but the resulting vector bundles on $C$ are stable for some 
polarisation.
\smallskip
ii)Assume $n'(n-n')(g-1)\le s$. Let $E\in U(n,d,(a_i))(C)$ be  generic .
 The set $A$ of subbundles of $E$ of rank 
$n'$ and degree $d'$ is non-empty of dimension $s-n'(n-n')(g-1)$ and the map
$$\matrix {\pi _Q:&A&\rightarrow &{\bf G}(n',E_Q)\cr
 &E'&\rightarrow& E'_Q\cr}$$
is either onto or has generic finite fibers.
                                                             
\bigskip
 
Proof: We assume $g\ge 2$. We shall deal with the case $g=1$ at the end.

1) Assume that $0<s\le n'(n-n')(g-1)$. Using (0,1),one can choose
 a vector bundle $E_g$
on $C_g$ of degree $d$ with $s_{n'}(E_g)=s$ and only a finite number of
vector bundles of rank $n'$ and maximal degree. Take as vector bundle 
$E_1$ on $C_1$ a direct sum of $n$ line bundles of degree zero.
Glue a direct sum of $n'$ of these line bundles with one of the maximal 
rank $n'$ subbundles of $E_g$. Then the resulting vector bundle
$E$ has a subbundle that gives the desired value $s$. It cannot have higher
degree subbundles because the highest degree for a vector bundle of $E_1$
is zero and we are assuming $s_{n'}(E_g)=s$. 
Moreover, the number of subbundles of maximum  degree is finite because
there are only a finite number of possibilities for the restriction of $E'$
to $C_1,C_g$. By (0,1), the conditions on stability are satisfied 
on $C_g$. By construction, they are satisfied on $C_1$.

2) Assume now $n'(n-n')(g-1)\le s\le n'(n-n')g, s\equiv n'd(n)$
. Take a direct sum of $n$ generic 
linebundles of degree zero on $C_1$. Take a generic vector bundle
$E_g$ on $C_g$ of degree $d$. Consider the set $A_g$ of subbundles of $E_g$
of rank $n'$ and degree $d'$.
  By (0,1), the image  $\pi _{g,Q}(A_g)$
 in ${\bf G}(n',E_{g,Q})$ has dimension $a_g(d')=s-n'n''(g-1)$.
 Consider the set $A_1$ of subbundles of $E_1$ of rank $n'$
and degree zero.  Then $A_1$ is finite and so is $\pi _1(A_1)$.
Consider the set $X$ of gluings of $E_{g,P_g}$ with $E_{1,P_1}$
that identify one direction in $\pi _{g,Q}(A_g)$ with one 
subspace in $\pi _1(A_1)$. By (3,1), the dimension of this set
is $a_g+0+n^2-n'(n-n')=s+n^2-n'(n-n')g$. Choose a generic element in 
$X$ to glue $E_{g,P_g}$ to $E_{1,P_1}$. The resulting $E$ has a 
subbundle of rank $n'$ and degree $d'$. By (0,1), the conditions 
on stability are satisfied on $C_g, g>1$. By (2,8) and (2,6), they are 
satisfied on $C_1$. By [T1] step 2, they give rise to a stable $E$.

 We need to show that $E$ 
has only a finite number of subbundles of rank $n'$ and degree
$d'$ and that it does not have subbundles of higher degree.

Assume it had a subbundle of higher degree. Denote by $\bar d_g$
and $\bar d_1$ the degrees of the restriction to $C_g$ and $C_1$
respectively. Then 
$\bar d_1\le 0.$
 Hence
 $$n'd-n\bar d_g=n'd-n(\bar d_g+\bar d_1)+n\bar d_1<s+n\bar d_1\le
 s\le n'(n-n')g.$$
Moreover, as $E_g$ is generic, $n'd-n\bar d_g \ge n'(n-n')(g-1)$(cf
[L] p.448 Satz). 
Hence, the induction assumption (in part ii)) of the Theorem applies
and the image $\pi_{P_g}(A(n',\bar d_g, E_g))$ in ${\bf G}(n',E_{P_g})$
 is $a_g(\bar d_g)=n'd-n\bar d_g-n'(n-n')(g-1)<n'(n-n')$.

As $n'd-n\bar d_g\ge n'(n-n')(g-1)$,
 $$n\bar d_1=n'd-n\bar d_g-
(n'd-n\bar d_g-n\bar d_1)>n'd-n\bar d_g-s\ge n'(n-n')(g-1)-s\ge -n'(n-n')$$
By (2,8) $\pi_{P_1}(A_1(n',\bar d_1,\bar E_1)$ 
has dimension $a_1=-n\bar d_1$. 
Now, $a_g(\bar d_g)+a_1(\bar d_1)=n'd-n\bar d_g-n\bar d_1-n'(n-n')(g-1)<
s-n'(n-n')(g-1)\le n'(n-n')$.
By (3,1), the dimension of the set of gluings that glue one subspace 
in $\pi _{g,Q}(A_g(n', \bar d_g, E_g))$ with one 
subspace in $\pi _1(\bar A_1(n',\bar d_1,E_1)$
has dimension at most 
$$a_g(\bar d_g)+a_1(\bar d_1)+n^2-n'(n-n')g=dimX.$$
Hence the generic element in $X$ does not glue any pair 
of these subspaces and no such subbundle of $E$ exists.

Assume now that $\bar d_g, \bar d_1$ are chosen in some other form so that 
$n'd-n\bar d_g-n\bar d_1-n'(n-n')(g-1)=s$. One can repeat the proof above
except that the strict inequalities are now equalities.
By (0,1), the fibers of $\pi_{g|A_g(n',\bar d_g, E_g)}$ 
are finite. By (2,8) , the fibers of $\pi_{1|A_1(n',\bar d_1, E_1)}$ 
are finite. By (3,1), a generic gluing in $X$ glues at most a finite number
 of subspaces in $\pi_g(A_g(n',\bar d_g, E_g))$ with subspaces in
$\pi_g(A_1(n',\bar d_1, E_1)).$ This shows that the number of
 subbundles obtained in this way is finite.

3)Assume now $n'(n-n')g\le s$.
Choose $ d_1, d'_1$ such that $n'd_1-nd'_1=n'(n-n')$.
This is possible because the greatest common divisor of 
$n,n'$ divides also $n'(n-n')$.

Write $h=(n,d_1), n=h \bar n, d_1=h\bar d_1$. Consider a vector 
bundle $E_1$ on $C_1$ that is the direct sum of $h$ indecomposable 
generic vector bundles of rank $\bar n$ and degree $\bar d_1$. 
Take a generic vector bundle $E_g$ on $C_g$ of rank $n$ and degree 
$d_g$. Glue $E_1,E_g$ at $P$ by a generic gluing. The resulting 
vector bundle $E$ is a generic point in one of the components of the
 moduli space of vector bundles of rank $n$ and degree $d$ on $C$
 (cf [T1] Theorem  ).

We first show that $A(n',d,E)$ is non-empty and has a component 
of dimension $s-n'(n-n')g$. Define $d'_g=d'-d'_1$. Then,
$n'd_g-nd'_g=s-n'(n-n')\ge n'(n-n')(g-1)$
By (0,1), $A_g(n',d_g,E_g)$ is non-empty and has 
a component of dim $n'd_g-nd'_g-n'(n-n')(g-1)$. The map to 
${\bf G}(n',E_{g,P_g})$ is either onto or it has finite fibers.
By (2,8), $A_1(n',d_1,E_1)$ is non-empty and has 
a component of dim $n'd_1-nd'_1=n'(n-n')$. The map to 
${\bf G}(n',E_{1,P_1})$ is onto and has finite fibers.
Therefore, each of the subbundles of degree $d'_g$ of $E_g$ glues with 
a finite number of subbundles of degree $d'_1$ of $E_1$. This 
shows that the component $A(n',(d'_1,d'_g),E)$ of $A(n',d,E)$ 
(corresponding to the given splitting of the degree of the subbundle 
between the  two components)is non-empty of 
dimension $n'd_g-nd'_g-n'(n-n')(g-1)=s-n'(n-n')g$.From (0,1), (2,8)
and (2,6), stability conditions are satisfied for the restrictions
 to the components. From [T1], they give rise to a stable $E$.

We need to check that $E$ does not have higher dimensional varieties of
 subundles of rank $n'$ and degree $d'$ corresponding to a different splitting
of the degree $d'=\bar d'_1+\bar d'_g$ .
Assume that $A(n',(\bar d'_1,\bar d'_g),E)$ is non-empty. By the genericity 
of $E_g$, $n'd_g-n\bar d'_g \ge n'(n-n')(g-1)$ (cf [L] ). 
By (0,1), $dimA(n',\bar d'_g,E_g)=n'd_g-n\bar d'_g-
n'(n-n')(g-1)$ and either is mapped by $\pi_{P_g}$ onto
 ${\bf G}(n',E_{g,P_g})$ or the map has finite fibers.
By (2,8), $dimA(n',\bar d'_1,E_1)=n'd_1-n\bar d'_1-
n'(n-n')(g-1)$ and either is mapped by $\pi_{P_1}$ onto
 ${\bf G}(n',E_{1,P_1})$ or the map has finite fibers.
An element in $A(n',(\bar d'_1,\bar d'_g),E)$ is determined 
by an element in $A(n',\bar d'_1,E_1)$ and an element in 
$A(n',\bar d'_g,E_g)$ whose images in ${\bf G}(n',E_{1,P_1})$
and ${\bf G}(n',E_{g,P_g})$ glue together. 

If $\pi_{P_g}$ is onto, then  each element in $A(n',\bar d_1,E_1)$ glues 
with a fiber of $\pi _{P_g}$. Hence we obtain a set of dimension
$n'd_1-n\bar d'_1+(n'd_g-n\bar d'_g-n'(n-n')=s-n'(n-n')g$ as desired.
If $\pi_{P_1}$ is onto, one can reason in the same way.
Assume now that both maps have finite fibers . Then by (3,1), the genericity of 
the gluing implies that $n'd_g-n\bar d'_g-n'(n-n')+n'd_1-n\bar d'_1+n^2
-n'(n-n')\ge n^2$. This can be written as $s\ge n'(n-n')g$.
Then, by (3,1)iv) the genericity of the gluing and the finiteness of the fibers 
of the maps $\pi_{P_1},\pi_{P_g}$, 
$$dimA(n',(\bar d_1,\bar d_g),E)=dim A_g(n',\bar d_g,E_g)+dim A_1(n',
\bar d_1,E_1)-n'(n-n')=s-n'(n-n')g.$$
 
Consider now the fiber $F$ of $\pi_{Q|A_1(n'\bar d_1, E_1)}$ over a
 generic point of ${\bf G}(n',E_{1,P_Q})$. From (2,8),
 either $\pi_{P_1}$ maps $F$
onto ${\bf G}(n',E_{1,P_1})$ or it maps with finite fibers.
As $A_g(n',\bar d_g, E_g)$ either maps onto ${\bf G}(n', E_{P_g})$ 
or it maps with finite fibers, the result follows with the same line
 of reasoning as before.

Consider now the case $g=1$. We can repeat the argument above
replacing (0,1) by (2,8) except that we have semistability instead 
of stability for the restriction of the vector bundle to $C_g=C_1'$.
 From the remark at the end of (1,2), we need to prove that 
the resulting vector bundle $E$ on the reducible curve is stable for some 
polarisation. If this were not the case, there would be a subbundle 
$F$ of $E$ with the same rank on each component and such that 
$\mu (F)=\mu (E)$(cf [T1], end of the proof of step 2, p.436.
 This latter condition along with the semistability
of each restriction implies that $\mu (F_{|C_1})=\mu (E_{|C_1}),
\mu (F_{|C_2})=\mu (E_{|C_2})$. There are only a finite number 
of subbundles $F'_1, F_1$ satisfying these conditions. 
From (3,1), the dimension of the set
of identifications of $E_{1,P_1}$ with $E'_{1,P'_1}$ 
that glue together a pair of these subspaces is at most $n^2-n'(n-n')$. 
This is smaller than the dimension of the $X$ considered above.
This finishes the proof. 
\bigskip 
The proof of (0,1) can now be carried out by induction  on $g$ 
using (1,2) and (3,2).
\bigskip

\beginsection References.

[A] M.Atiyah, Vector bundles on elliptic curves, Proc.London Math.Soc.
{\bf (3)7}(1957),414-452.

\noindent [B,B,R] E.Ballico, L.Brambila-Paz, B.Russo. Exact sequences 
of stable vector bundles on projective curves, Math.Nachr.

\noindent [L]H.Lange, Zur Klassifikation von Regelmanigfaltigkeiten, Math.Ann.
{\bf 262}(1983), 447-459.

\noindent [L2] H.Lange, Universal Families of Extensions, 
J. of Algebra {\bf 83} (1983), 101-112.

\noindent [L,N] H.Lange, M.S.Narasimhan, Maximal subbundles of rank
two vector bundles on curves, Math.Ann. {\bf 266} (1983), 55-72.

\noindent [P] R.Pandharipande, A compactification over $\bar {\cal M}
_g$ of the universal moduli space of slope-semistable vector 
bundles, Journal of the A.M.S. {\bf V9 N2}(1996), 425-471.

\noindent [S] C.S.Seshadri. Fibres vectoriels sur les courbes 
algebriques. Asterisque {\bf 96}, Societe Mathematique de France, 
1982.

\noindent [T1]M. Teixidor i Bigas, Moduli spaces of semistable bundles on 
tree-like curves, Math.Ann.{\bf 290}(1991), 341-348.

\noindent [T2]M. Teixidor i Bigas, Moduli spaces of semistable bundles on 
reducible curves, Amer.J.of Math. {\bf 117} (1995), 125-139

\noindent [T3]M. Teixidor i Bigas, Stable extensions by line bundles,
preprint Jan. 97.

\noindent [Tu]L.Tu, Semistable bundles over an elliptic curve, Adv in Math
{\bf 98}(1993),1-26.

\end